\documentclass[12pt,reqno]{amsart}

\usepackage[noadjust]{cite}

\usepackage{hyperref}
\usepackage{orcidlink}

\definecolor{sinkgreen}{RGB}{220,247,220}
\definecolor{sourcered}{RGB}{250,220,220}

\tikzset{
  simple vertex/.style={
    circle,
    draw,
    thick,
    fill=white,
    minimum size=9mm,
    inner sep=1pt
  },
  path vertex/.style={
    rectangle,
    rounded corners=2pt,
    draw,
    thick,
    fill=white,
    minimum width=4.1cm,
    minimum height=10mm,
    align=center,
    inner sep=3pt
  },
  sink/.style={fill=sinkgreen,very thick},
  source/.style={fill=sourcered,very thick},
  fone/.style={thick},
  ftwo/.style={thick,dashed},
  oriented/.style={-{Stealth[length=2.4mm,width=1.8mm]}},
  edge label/.style={fill=white,inner sep=1.5pt,font=\small}
}

\usepackage{amsmath,amssymb}
\usepackage{xcolor}
\usepackage{tikz}
\usetikzlibrary{arrows.meta,positioning}

\usepackage{xcolor}

\usepackage{amsmath,amsthm,amsfonts,amssymb,amscd,  multicol, verbatim, enumerate, stackrel}
\usepackage{verbatim}
\usepackage{fullpage}
\input{xy}
\usepackage{xcolor}
\usepackage[english]{babel}
\usepackage{cite}
\usepackage{amsfonts}
\usepackage{latexsym}
\usepackage{amsthm}
\usepackage{amsopn}
\usepackage{amsmath}
\usepackage[all]{xy}

\usepackage{hyperref}

\usepackage{verbatim}
\usepackage{amssymb}
\usepackage{mathrsfs}
\usepackage{float}
\usepackage{fullpage}
\usepackage{amscd}
\usepackage{amscd}
\usepackage{color}

\usepackage{graphicx}
\usepackage{tikz}

\theoremstyle{plain}
\newtheorem*{thm*}{Theorem}
\newtheorem*{con*}{Conjecture}
\newtheorem{thm}{Theorem}[section]
\newtheorem{pro}[thm]{Proposition}

\newtheorem{lem}[thm]{Lemma}
\newtheorem{con}[thm]{Conjecture}

\newtheorem{cor}[thm]{Corollary}
\newtheorem{exa}[thm]{Example}

\newtheorem*{str*}{Strategy}

\theoremstyle{definition}
\newtheorem{defn}[thm]{Definition}

\theoremstyle{remark}
\newtheorem{rmk}[thm]{Remark}

\usepackage[noEnd=true,indLines=false]{algpseudocodex}

\newcommand{\D}{\mathbb{D}}
\title[]{Flip Combinatorial Invariance for flipclasses with small residues}
\author[Esposito, Marietti, Paolini]{Francesco Esposito\, \orcidlink{0000-0003-4627-9193}, Mario Marietti\, \orcidlink{0000-0001-7045-0347}, Giovanni Paolini\, \orcidlink{0000-0002-3964-9101}}

\date{}

\address{Francesco Esposito, Dipartimento di Matematica, Universit\`a degli Studi di Padova,
via Trieste 63, 35121 Padova, Italy}
\address{Mario Marietti, Dipartimento  di Ingegneria Industriale e Scienze Matematiche, Universit\`a Politecnica delle Marche, via Brecce Bianche, 60131 Ancona,  Italy}
\address{Giovanni Paolini, Dipartimento di Matematica, Universit\`a di Bologna, Piazza di Porta San Donato 5, 40126 Bologna, Italy}

\email{esposito@math.unipd.it}
\email{m.marietti@univpm.it}
\email{g.paolini@unibo.it}


\begin{document}

\maketitle

\begin{abstract}
We establish (Normal) Flip Combinatorial Invariance for flipclasses with no residues longer than 8. Actually, we show more generally that all flipclasses of finite type with no residues longer than 8 have only one increasing path. For flipclasses of Weyl type, these are precisely the flipclasses having only one increasing path. This marks a refinement of the classical Combinatorial Invariance of Kazhdan--Lusztig polynomials.
\end{abstract}

\section{Introduction}
The  Combinatorial Invariance Conjecture,   independently formulated by Lusztig and Dyer in the 1980s,  predicts that the Kazhdan--Lusztig polynomial associated with a Bruhat interval $I$ in a Coxeter group depends only on the isomorphism class (as a poset) of $I$. 
This conjecture, which has been extensively studied over the past decades (\cite{B-G,BGL,BBDVW,Bre03,Bre04,BCM1,BCM2, B-M,BLP,CM, CMult,Del,dC,EM,EM2,EMS,G-W,Inc,Inc2,
MJoA,MJaco,Mtrans,M,P}), has led to the introduction of new interesting combinatorial concepts and new conjectures implying the Combinatorial Invariance Conjecture.

One of these conjectures, the Flip Combinatorial Invariance Conjecture \cite{EMS}, stems from the idea that Bruhat intervals should be decomposed into  simpler constituents, namely the flipclasses. In \cite{EMS},
for each $u$ and $v$  in an arbitrary Coxeter group $W$, with $u\leq v$ in Bruhat order, and each $h \in \mathbb N$,  the set of paths (in the Bruhat graph of $W$) from $u$ to $v$ of length $h$ is partitioned  into smaller sets, the {\em $h$-flipclasses} from $u$ to $v$, and a notion of combinatorial isomorphism of flipclasses is introduced.
The number of increasing paths of a flipclass $F$ of finite type (so, in particular, of Weyl type) is shown to be independent of the chosen reflection ordering in \cite[Theorem~5.1/Proposition 6.7]{EMS}. It is denoted $c(F)$.
\begin{con}[Flip Combinatorial Invariance Conjecture]
If two flipclasses $F$ and $F'$ of Weyl type are combinatorially isomorphic, then $c(F)=c(F')$. 
\end{con}

In this work, we introduce a new edge-labelled graph attached with a flipclass $F$, the \emph{flipgraph} of $F$, denoted $FG_F$. If two   flipclasses are combinatorially isomorphic, then their flipgraphs are isomorphic. 
The following conjecture is stronger than the Flip Combinatorial Invariance Conjecture, and hence stronger than the Combinatorial Invariance Conjecture (for Weyl groups). 

\begin{con}[Normal Flip Combinatorial Invariance Conjecture]
If two flipclasses $F$ and $F'$ of Weyl type  have isomorphic flipgraphs, then $c(F)=c(F')$. 
\end{con}

\bigskip
We provide a simple characterization of the flipclasses of finite type having exactly one increasing path. This characterization is given in terms of certain elementary subgraphs of the flipgraph $FG_F$, viz. its \emph{residues}, which are the cycles of $FG_F$ obtained by following only two labels.  For flipclasses of Weyl type, the characterization is combinatorial and in fact depends only on the flipgraph. Hence, it establishes Normal Flip Combinatorial Invariance for flipclasses of Weyl type having exactly one increasing path.

In the following theorem, we collect the main results of this work.
\begin{thm}
\label{collect}
$ $
\begin{enumerate}
\item Let $F$ be a flipclass of finite type. Then $F$ has only one increasing path if and only if each residue $S$ of $F$ satisfies $c(S)=1$. In particular, if each residue of $F$ is a $2m$-cycle with $m\leq 4$, then $c(F)=1$.
\item Let $F$ be a flipclass of Weyl type. Then $F$ has only one increasing path if and only if each residue $S$ of $F$ is a $2m$-cycle with $m\leq 4$.
\end{enumerate}
\end{thm}

We state explicitly the following consequence, which establishes a first highly nontrivial instance of Normal Flip Combinatorial Invariance. 
\begin{cor}
Let $F$ and $F'$ be two flipclasses of Weyl type with isomorphic flipgraphs. Then 
$$c(F)=1 \iff c(F')=1.$$
\end{cor}

\bigskip

In October 2026, OpenAI released a
large collection of mathematical manuscripts produced
by an internal frontier model.
One of these manuscripts (see \cite{O}) claims to prove
the Combinatorial Invariance Conjecture.
%asserts that Kazhdan--Lusztig polynomials are determined by the poset structure of Bruhat intervals
The Flip Combinatorial Invariance Conjecture remains unsettled.

\bigskip
\noindent {\bf Disclosure of AI use}. 
In the preparation of this work, we benefited from interaction with GPT-6 Astra (OpenAI) and Opus 5.5 (Anthropic) working in parallel and coordinating themselves. More precisely, they suggested what follows.
\begin{itemize}
\item A sketch of the proof of the fact that, when a flipclass of  type $A$ has no residues of length 10, then $c(F)=1$. This result is not included in the paper because it is weaker than Theorem~\ref{minore2}, (\ref{minoreft}) but the sketched proof was useful to find the proof of the formula of Theorem~\ref{minore2}, (\ref{minoreft}), which was also suggested (without proof).
\item Theorem~\ref{terni}, both its statement and its proof.
\end{itemize}
The code for the final check of  Proposition~\ref{2,3,4,5,6,7,9} and the TikZ codes used to produce the figures were generated with the assistance of Opus~5.5 and subsequently checked by the authors. That agent was also used for spell-checking.
The authors take full responsibility for the correctness of the results and the quality of their presentation.

\section{Basics}

In this section, we  review some background material. 

For $n\in \mathbb N^+$, we let $[n]=\{1,2,\ldots,n\}$.

\subsection{Coxeter groups}
We refer the reader to  \cite{BB} for the basic facts regarding Coxeter group theory.

The \emph{Bruhat graph} of $W$, denoted $B(W)$, is the edge-labelled directed graph whose vertex set is $W$ and where  $u \xrightarrow{t} v$ if and only if $\ell(u)<\ell(v)$ and $vu^{-1}=t \in T$. The \emph{Bruhat order} is the partial order on $W$ where $u \leq v$ whenever there is a (directed) path from $u$ to $v$ in $B(W)$. 
 Let $u,v \in W$, $u \leq v$, and  $\Gamma= (u=x_0 \stackrel{t_1}{\longrightarrow} x_1 \stackrel{t_2}{\longrightarrow}  \cdots \stackrel{t_{h-1}}{\longrightarrow} x_{h-1}\stackrel{t_h}{\longrightarrow}  x_h=v)$ be a path from $u$ to $v$ in $B(W)$. We define
  $h$ to be the {\em length} of $\Gamma$. 
 For short, we say  $\Gamma$ passes through $x_i$ at time $i$.

Let $W'$ be any reflection subgroup of a Coxeter group $W$, i.e., $W'= \langle W' \cap T\rangle$. Then $W'$ is itself a Coxeter group with the set $\{ t \in T : \{r\in T \cap W' : \ell(rt)<\ell(t) \}=\{t\}  \}$, denoted $\chi(W')$, as a set of Coxeter generators  (see \cite{DyeJAlg}). 

Recall (see \cite[\S 2]{DyeComp}) that a {\em reflection ordering} on $(W,S)$ is a total order $\preceq$ on $T$ such that if $(W',\chi(W'))$
is a dihedral reflection subgroup of $W$ (i.e. $\chi(W')$
has cardinality $2$) then either
\[
a \preceq aba \preceq ababa 	\preceq \cdots \preceq babab \preceq bab \preceq b
\]
or 
\[
b \preceq bab \preceq babab \preceq \cdots \preceq ababa \preceq aba \preceq a
\]
where $\{ a,b \}= \chi(W')$ and the expressions appearing above are reduced and exhaust the reflections in $W'$. 
Given a reflection ordering $\preceq $,  a  path $\Gamma= \big(x_0 \stackrel{t_1}{\longrightarrow} x_1 \stackrel{t_2}{\longrightarrow}  \cdots \stackrel{t_{h-1}}{\longrightarrow} x_{h-1}\stackrel{t_h}{\longrightarrow} x_h\big)$ in $B(W)$ is \emph{increasing with respect to $\preceq $} provided that $t_1 \preceq t_2\preceq \cdots \preceq t_h$.

\subsection{Flipclasses for arbitrary Coxeter groups}
Let us recall some definitions from \cite{EMS}. 
Let $W$ be an  arbitrary Coxeter group.  
Given $h \in \mathbb N$ and $u,v \in  W$, we denote the set of $h$-paths from $u$ to $v$ in the Bruhat graph $B(W)$ by $P_{h}(u,v)$. 

\begin{lem}[Lemma~3.1 of \cite{EMS}]
\label{mari}
Let $u,v \in W$. 
\begin{enumerate}
\item $P_2(u,v)$ contains an even number of paths.
\item Let $X=\{x\in [u,v]\setminus\{u,v\} : \text{ there exists $\Gamma$ in $P_2(u,v)$ containing $x$} \}$.  Then there exists a unique fixed-point free involution  $f :X \to X$ such that, given $x_1,x_2 \in X$, with $\ell(x_1) < \ell(x_2)$ and $x_2\neq f(x_1)$, we have $\ell(f(x_1)) < \ell(f(x_2))$.  
\end{enumerate} 
\end{lem}
\begin{defn}
\label{flipflip}
The \emph{flip} of a path $ (u{\longrightarrow} x {\longrightarrow} v)$ is the path  $(u{\longrightarrow} f(x) {\longrightarrow} v)$, where $f$ is the involution of Lemma~\ref{mari}.  

Let $u,v \in W$  and $h \in \mathbb N$. 
\begin{itemize}
\item
The {\em $i$-th flip} is the map $f_i: P_{h}(u,v) \to P_{h}(u,v)$ that sends a path $\Gamma = \big(u=x_0 {\longrightarrow} x_1 {\longrightarrow} \cdots {\longrightarrow} x_{h}=v\big)$ to the path obtained from $\Gamma$ by substituting the subpath $( x_{i-1} {\longrightarrow} x_i {\longrightarrow}  x_{i+1})$ with its flip, for $i\in [h-1]$.
\item 
An \emph{$h$-flipclass} of $W$ from $u$ to $v$ is a minimal subset of $P_{h}(u,v)$ closed under flips.
\end{itemize}
\end{defn}

The definition of $i$-th flip, and hence the definition of a flipclass, depends only on the isomorphism type of the Bruhat graph.

Given a flipclass $F$, we let $E(F) =\{w\in W : \textrm{ there exists $\Gamma$ in $F$ with $w\in \Gamma$}\}$. 
\begin{defn}
Let $W$ and $W'$ be two Coxeter groups,  and $h \in \mathbb N$. Let $F$ be an $h$-flipclass  of $W$ and $F'$ be an $h$-flipclass  of $W'$. 
 We say that $F$ and $F'$ are \emph{combinatorially isomorphic}  provided that there exists a bijection  $\varepsilon$  from $E(F)$ to $E(F')$ such that:
%\begin{itemize}   
 $$\big(x_0{\longrightarrow} x_1 {\longrightarrow}  \cdots {\longrightarrow} x_{h-1}{\longrightarrow} x_h\big) \in F$$ if and only if  
$$\big(\varepsilon(x_0) {\longrightarrow} \varepsilon(x_1) {\longrightarrow}  \cdots {\longrightarrow} \varepsilon(x_{h-1}){\longrightarrow} \varepsilon(x_h)\big) \in F'$$ 
and  the induced  bijection $\bar{\varepsilon}$ from $F$ to $F'$ is flip-preserving, i.e., $\bar{\varepsilon}f_i(\Gamma)= f_{i} \bar{\varepsilon}(\Gamma)$  for all $\Gamma \in F$.
\end{defn}

\begin{rmk}
\label{costante}
Let $[u,v]$ and $[u',v']$ be two isomorphic Bruhat intervals in two (possibly different) Coxeter groups. The multiset of combinatorial isomorphism classes of the $h$-flipclasses from $u$ to $v$ and that  of combinatorial isomorphism classes of the $h$-flipclasses from $u'$ to $v'$ coincide.
\end{rmk}

\begin{defn}
\label{deffin}
Let $F$ be a flipclass of an arbitrary Coxeter group and let $W(F)$ be the reflection subgroup generated by the reflections labeling the paths in $F$. Then we say that $F$ is of a certain type if $W(F)$ is a Coxeter group of that type.
%\begin{itemize}
%\item \emph{finite type},  if $W(F)$ is a finite Coxeter group;
%\item \emph{Weyl type}, if $W(F)$ is  a Weyl group;
%\item  \emph{symmetric type}, if $W(F)$ is a product of Weyl  groups of type $A$;
%\item \emph{dihedral type},  if $W(F)$ is a dihedral Coxeter group.
%\end{itemize}
\end{defn}
Clearly, a flipclass of type $A$ is of Weyl type, and a flipclass of Weyl type is of finite type.

We will need the following lemma, where, for two sets $A$ and $B$, the notation $A\preceq B$ is shorthand for $a\preceq b$ for all $a\in A$ and $b\in B$.
\begin{lem}[Lemma 3.20 of \cite{EMS}]
\label{espomariancona24}
Let $u,v \in W$ and let $\preceq$ be a reflection ordering of $W$. Let $(u  \stackrel{t_1} {\longrightarrow} x  \stackrel{t_2}{\longrightarrow} v)$ and $(u \stackrel{r_1}{\longrightarrow} y  \stackrel{r_2}{\longrightarrow} v)$ be one the flip of the other. Then one of the two  paths is increasing and the other is decreasing. Moreover, either $\{t_1,r_2\} \preceq \{t_2,r_1\}$, or  $\{t_2,r_1\} \preceq \{t_1,r_2\}$. 
\end{lem}

Given a reflection ordering $\preceq$  of $W$ and a flipclass $F$, we denote by $c_{\preceq}(F)$ the number of paths of $F$ that are increasing with respect to $\preceq$. By the following result, $c_{\preceq}(F)$ is always at least 1.
\begin{pro}[Proposition 3.21 of \cite{EMS}]
\label{greedy}
Fix a reflection ordering   $\preceq$ of $W$. Let $F$ be a flipclass. The  lexicographically first path $\Gamma$ of $F$ is increasing. 
\end{pro}

We end this section with the following result. 
\begin{thm}[Theorem 5.1/Proposition 6.7 of \cite{EMS}]
\label{nondipende}
Let $F$ be a  flipclass of finite type. Then the  number $c_{\preceq}(F)$ does not depend on the reflection ordering.
\end{thm}
After Theorem \ref{nondipende}, given a flipclass of finite type,  we denote by $c(F)$ the number of increasing paths with respect to any reflection ordering.

\section{Subflipclasses, flipgraphs, and residues}
In this section, we introduce some concepts that are central in this work.

Let $W$ be an arbitrary Coxeter group.
\subsection{Subflipclasses}
In this subsection, we introduce the notion of a subflipclass of an $h$-flipclass $F$ of $W$.
\begin{defn}
For $J \subseteq [h-1]$, a \emph{$J$-subflipclass} of $F$ is  a subset $O$ of $F$ that is an orbit of
$\langle f_j : j \in J\rangle$.
\end{defn} 
Notice that the paths of a $J$-subflipclass $O$ agree at every time not in $J$.  

If $J$ is an interval, we call $O$ a \emph{segment subflipclass}. In this case, we let $\bar{O}$ be the $(|J|+1)$-flipclass given by forgetting the common initial and the common final parts. Clearly, if $F$ is of finite type (respectively, Weyl type or type $A$), then $\bar{O}$ is of finite type (respectively, Weyl type or type $A$).

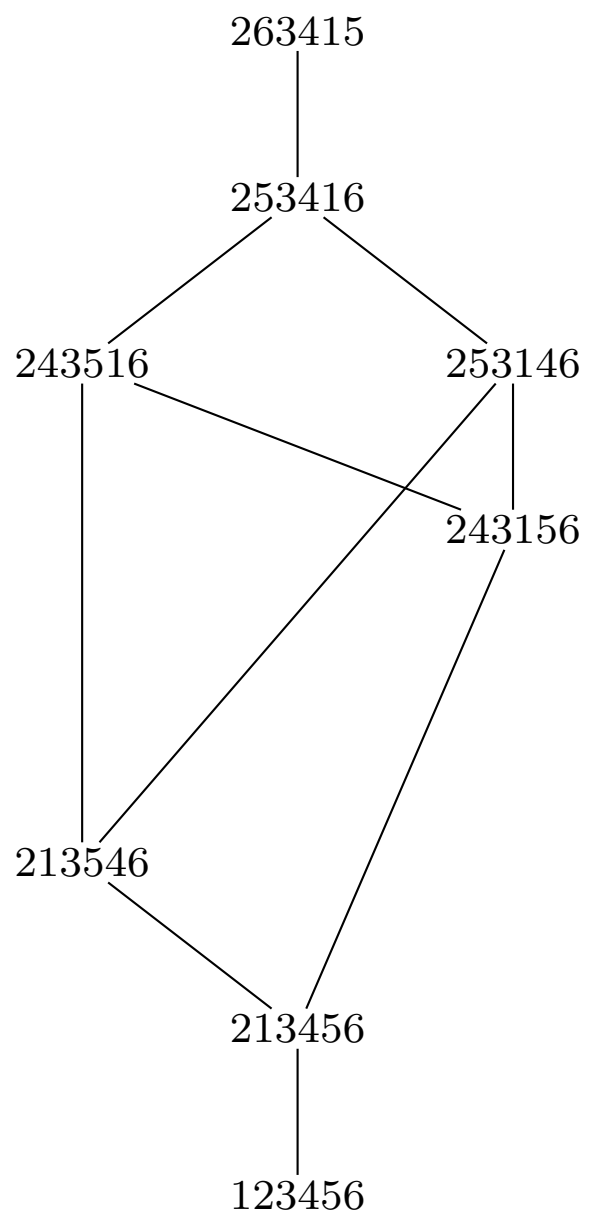
\begin{figure}[ht]
\centering

\begin{tikzpicture}[every node/.style={font=\tiny, inner sep=1pt}]
  \node (n123456) at (0.000,0.000) {123456};
  \node (n213456) at (0.000,1.100) {213456};
  \node (n213546) at (-1.425,2.200) {213546};
  \node (n243156) at (1.425,4.400) {243156};
  \node (n243516) at (-1.425,5.500) {243516};
  \node (n253146) at (1.425,5.500) {253146};
  \node (n253416) at (0.000,6.600) {253416};
  \node (n263415) at (0.000,7.700) {263415};
  \draw (n123456) -- (n213456);
  \draw (n213456) -- (n213546);
  \draw (n213456) -- (n243156);
  \draw (n213546) -- (n243516);
  \draw (n213546) -- (n253146);
  \draw (n243156) -- (n243516);
  \draw (n243156) -- (n253146);
  \draw (n243516) -- (n253416);
  \draw (n253146) -- (n253416);
  \draw (n253416) -- (n263415);
\end{tikzpicture}
\caption{A segment subflipclass}
\label{notte}
\end{figure}

\begin{exa}
Let $F$ be the $5$-flipclass of the path $\Gamma = 123456 \longrightarrow 213456 \longrightarrow 213546 \longrightarrow 243516 \longrightarrow 253416 \longrightarrow 263415$ in $S_6$. Let $J = \{2,3\}$. Then the $J$-subflipclass $O$ of $F$ containing $\Gamma$ has the four paths shown in Figure~\ref{notte}. The subflipclass $O$ is a segment subflipclass, and $\overline{O}$ is a 3-flipclass of dihedral type.
\end{exa}

Let $O$ be a $J$-subflipclass. If $I_1,\ldots,I_r$ are the maximal intervals of $J$, then $O$ may be identified with 
the product $\bar{O}_1 \times \cdots \times \bar{O}_r$.
For a reflection ordering $\preceq$, we let 
$$c_{\preceq}(O)= \prod_{i=1}^{r} c_{\preceq}(\bar{O}_i).$$

Notice that one has $c_{\preceq}(O) \geq 1$, by Proposition~\ref{greedy}. When we are in finite type, we drop the subscript $\preceq$.
%An isomorphism of flip graphs (hence, a fortiori, a  combinatorial isomorphism of flipclasses) maps subflipclasses to subflipclasses and segment classes to isomorphic segment classes.

\begin{defn}
We say that an $h$-flipclass $F$ is \emph{normal} if, for each subset $J$ of $[h-1]$, given two paths $\Gamma$ and $\Gamma'$ of $F$ that pass through the same element at time $i$ for every $i\in [h]\setminus J$, then $\Gamma$ and $\Gamma'$ are contained in the same $J$-subflipclass.
\end{defn}

\subsection{The flipgraph of a flipclass}
In this subsection, we introduce the flipgraph of a flipclass.

\begin{defn}
 Let $F$ be an $h$-flipclass. The \emph{flipgraph} of $F$, denoted $FG_F$, is the graph with $F$ as vertex set and having an edge labeled $i$ between the paths $\Gamma$ and $f_i(\Gamma)$, for each $\Gamma \in F$ and $i\in[h-1]$.
\end{defn}
\begin{defn}
For every choice of a reflection ordering $\preceq$, we can direct the edges of $FG_F$ from the path with a descent at the flip index to the path with an ascent there. This is well defined by Lemma~\ref{espomariancona24}. We denote by $_{\preceq} FG_F$ the resulting directed edge-labeled graph, and we call it the \emph{$\preceq$-oriented flipgraph} of $F$.
\end{defn}

An isomorphism of flipgraphs (respectively, oriented flipgraphs) is an isomorphism of graphs (respectively, oriented graphs) preserving edge labels. 
The isomorphism class of the undirected graph $FG_F$ depends only on the combinatorial isomorphism class of $F$. 
The isomorphism class of the directed graph $_{\preceq} FG_F$ depends on the choice of the reflection ordering (see Example~\ref{diversi}).

\begin{lem}\label{loreto}
 Let $F$ be a flipclass and let $\preceq$ be a reflection ordering.
\begin{enumerate}
\item The $\preceq$-increasing (respectively,  $\preceq$-decreasing) paths are precisely the sinks (respectively, the sources) of  $_{\preceq} FG_F$.

\item The oriented graph $_{\preceq} FG_F$ has no (oriented) cycles; it is thus endowed with a partial order.

\item Maximal oriented paths of $_{\preceq} FG_F$ start at sources and end at sinks.
\end{enumerate}
\end{lem}
\begin{proof}
The first statement is clear from the definitions.
By Lemma~\ref{espomariancona24}, the source of an oriented edge is lexicographically greater than the target of that same edge. Thus, there are no oriented cycles and the transitive closure of the relation given by the oriented edges (i.e., $u$ is related to $v$ if and only if there is a directed path from $u$ to $v$) is a partial order refined by the lexicographical order.
This proves the second statement.
The third statement readily follows from the second statement and the definitions.
\end{proof}

\begin{figure}[ht]
\centering
% #1 = horizontal shift, #2 = vertices (index/style/descent set),
% #3 = edges (from/to/label/style), listed clockwise starting from P_0--P_1,
% #4 = name of the graph
\newcommand{\flippanel}[4]{%
  \begin{scope}[xshift=#1]
    \foreach \k/\s/\D in {#2} {
      \pgfmathsetmacro{\ang}{90-60*\k}
      \node[\s, label={[descent]\ang:$\D$}] (v\k) at (\ang:1.3) {$P_\k$};
    }
    \foreach \a/\b/\i/\s [count=\e] in {#3} {
      \pgfmathsetmacro{\ang}{120-60*\e}
      \draw[\s] (v\a) -- (v\b);
      \node[edgelabel] at (\ang:0.83) {$\i$};
    }
    \node at (0,-2.4) {$#4$};
  \end{scope}%
}
\begin{tikzpicture}[
  vertex/.style={circle, draw, thick, fill=white, minimum size=6.5mm,
                 inner sep=0pt, font=\footnotesize},
 % sink/.style={vertex, draw=green!45!black, fill=blue!30},
 sink/.style={vertex, draw=blue!45!black, fill=blue!40},
  source/.style={vertex, draw=red!60!black, fill=red!40},
  plain/.style={thick},
  arrow/.style={thick, -{Stealth[length=2.1mm,width=1.6mm]}},
  edgelabel/.style={inner sep=0pt, font=\scriptsize},
  descent/.style={font=\scriptsize, inner sep=1.5pt, text=black!70},
]
\flippanel{0cm}
  {0/vertex/{}, 1/vertex/{}, 2/vertex/{}, 3/vertex/{}, 4/vertex/{}, 5/vertex/{}}
  {0/1/1/plain, 1/2/2/plain, 2/3/1/plain, 3/4/2/plain, 4/5/1/plain, 5/0/2/plain}
  {FG_F}
\flippanel{4.2cm}
  {0/sink/{\varnothing}, 1/vertex/{\{1\}}, 2/vertex/{\{2\}},
   3/source/{\{1,2\}}, 4/vertex/{\{1\}}, 5/vertex/{\{2\}}}
  {1/0/1/arrow, 2/1/2/arrow, 3/2/1/arrow, 3/4/2/arrow, 4/5/1/arrow, 5/0/2/arrow}
  {{}_{\preceq}FG_F}
\flippanel{8.8cm}
  {0/sink/{\varnothing}, 1/vertex/{\{1\}}, 2/source/{\{1,2\}},
   3/vertex/{\{2\}}, 4/vertex/{\{1\}}, 5/vertex/{\{2\}}}
  {1/0/1/arrow, 2/1/2/arrow, 2/3/1/arrow, 3/4/2/arrow, 4/5/1/arrow, 5/0/2/arrow}
  {{}_{\preceq'}FG_F}
\end{tikzpicture}
\caption{The flipgraph $FG_F$ and its orientations with respect to the
reflection orderings $\preceq$ and $\preceq'$. Sinks are shown in a hue of  blue and sources
in a hue of red; the set next to each vertex is its descent set.}
\label{firenze}
\end{figure}
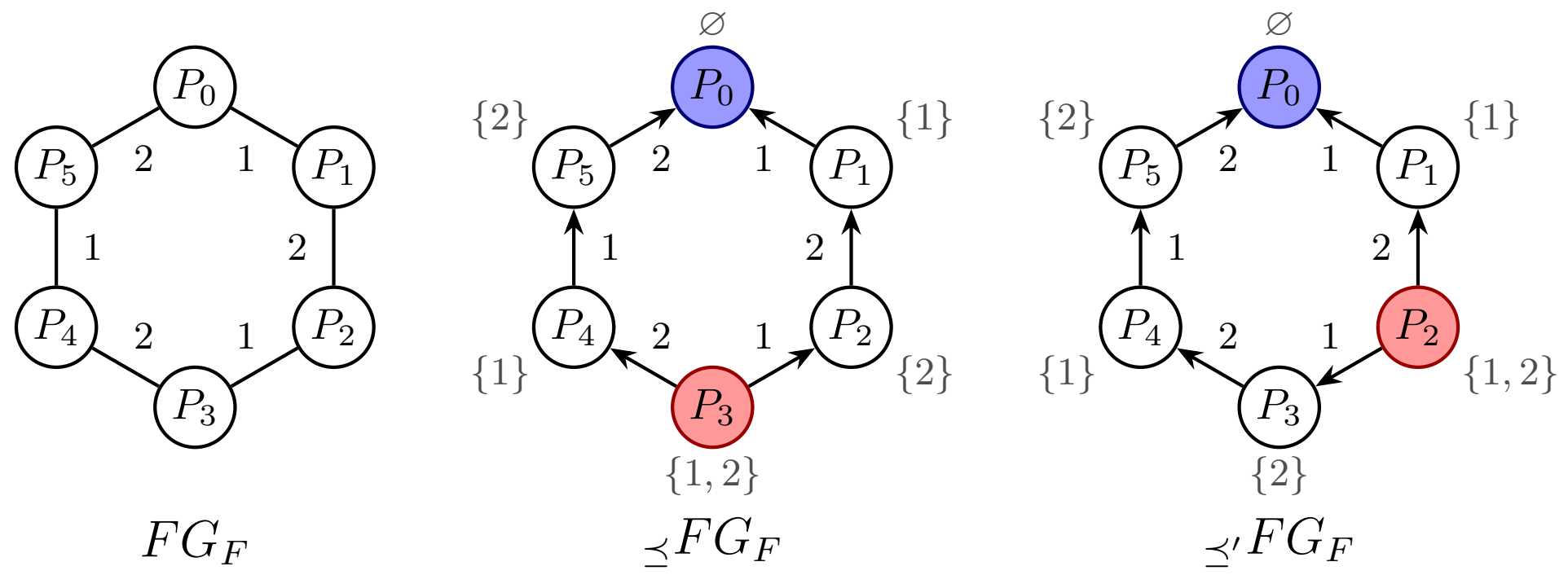

\begin{exa}\label{diversi}
\newcommand{\flipto}[1]{\xrightarrow{(#1)}}%
Let $F$ be the flipclass of $S_4$ consisting of the following six paths:
\[
\begin{array}{@{}r@{\;=\;}l@{}}
P_0 & 1234 \flipto{12} 2134 \flipto{13} 2314 \flipto{14} 2341,\\
P_1 & 1234 \flipto{23} 1324 \flipto{12} 2314 \flipto{14} 2341,\\
P_2 & 1234 \flipto{23} 1324 \flipto{24} 1342 \flipto{12} 2341,\\
P_3 & 1234 \flipto{34} 1243 \flipto{23} 1342 \flipto{12} 2341,\\
P_4 & 1234 \flipto{34} 1243 \flipto{12} 2143 \flipto{13} 2341,\\
P_5 & 1234 \flipto{12} 2134 \flipto{34} 2143 \flipto{13} 2341.
\end{array}
\]
The flipgraph $FG_F$ of $F$ is the graph on the left in Figure~\ref{firenze}.
Consider the two reflection orderings
\begin{align*}
\preceq\,  &:\quad (12)\preceq(13)\preceq(23)\preceq(14)\preceq(24)\preceq(34),\\
\preceq'\, &:\quad (12)\preceq'(13)\preceq'(14)\preceq'(34)\preceq'(24)\preceq'(23).
\end{align*}
The oriented flipgraphs ${}_{\preceq}FG_F$ and ${}_{\preceq'}FG_F$ are the two
graphs on the right in Figure~\ref{firenze}, and they are not isomorphic.
\end{exa}

While the isomorphism class of the directed graph $_{\preceq} FG_F$ depends on the choice of the reflection ordering, the number of sinks and the number of sources (which coincide) do not for flipclasses of finite type (see \cite[Theorem~5.1/Proposition~6.7]{EMS}).

Informally, if there is a directed path from $u$ to $v$, then we say that $v$ is downstream from $u$, and $u$ is upstream from $v$. 

\begin{defn}
\label{omig}
Given a vertex $v$, we let $\Omega(v)$ be the set of sinks that are downstream from $v$. 
\end{defn}
Note that if $u$ is upstream from $v$, then $\Omega(u)\supseteq \Omega(v)$.

Fix a reflection ordering $\preceq$.  We let $\mathcal S^{\preceq}_{u}$ (respectively, $\mathcal S^{\preceq}_{d}$) be the set of sources (respectively, sinks) of $_{\preceq} FG_F$. We drop the dependence from $\preceq$ in the notation when no confusion arises. Consider the following bipartite graph. The vertex set is $\mathcal S_{u} \cup \mathcal S_{d}$. Each  edge connects a source and a sink. A source $u$ and a sink $v$ are adjacent  if and only if $v\in \Omega(u)$.
\begin{lem}
\label{bipartito-connesso}
The above graph is connected. 
\end{lem}
\begin{proof}
Let $D_1 \sqcup D_2= \mathcal S_{d}$ be a non trivial partition of $ \mathcal S_{d}$. In order to prove the assertion, we show that there exists $s$ in $\mathcal S_{u}$ such that $\Omega(s)\cap D_i\neq \emptyset$, for $i=1,2$. 

Consider $U_1= \{v \in  _{\preceq} FG_F: \Omega(v) \cap D_1 \neq \emptyset\}$ and   $U_2=\{v \in  _{\preceq} FG_F: \Omega(v) \cap D_2 \neq \emptyset\}$. By Lemma~\ref{loreto}, we have that $U_1 \cup U_2$ is the vertex set of $_{\preceq} FG_F $. 
We claim $U_1 \cap U_2 \neq \emptyset$. Indeed, if $U_1 \cap U_2 = \emptyset$, we would have a partition of the vertex set. 
Since $ FG_F$ is connected, we would have an arc from a vertex $v$ in one of the two sets to a vertex in the other set; such vertex $v$ would belong to both $U_1$ and $U_2$. So the claim is proved.

Let $v\in U_1 \cap U_2$. By Lemma~\ref{loreto}, there exists a source $s$ upstream from $v$, and thus  $\Omega(s)\cap D_i\neq \emptyset$, for $i=1,2$. Hence, we have the assertion. 
\end{proof}

\subsection{Residues of a flipclass}
We now introduce the residues of a flipclass $F$, which actually depend only on the flipgraph of $F$.
\begin{defn}
Let $F$ be an $h$-flipclass. An \emph{$\{i,j\}$-residue} of  $F$ (or also of $FG_F$) is a connected component of the labeled graph that one obtains from $FG_F$ by deleting each edge with label different from $i$ and $j$.
\end{defn}
We say that an $\{i,j\}$-residue is \emph{adjacent} provided that $|i - j| = 1$, \emph{non-adjacent} otherwise. Note that each residue is a cycle of even length. 

\begin{pro}
\label{corrisponde}
Let $F$ be a flipclass.
\begin{enumerate}
\item Every non-adjacent residue $S$ of $F$ is a 4-cycle.
\item Every adjacent residue $S$ of $F$ corresponds to a segment subflipclass $O$ of $F$ such that $\bar{O}$ is a $3$-flipclass.
\end{enumerate}
\end{pro}
\begin{proof}
Straightforward.
\end{proof}

Let $F$ be a flipclass and $S$  be one of its  residues. Let $\preceq$ be a reflection ordering. We denote by $_{\preceq}S$ the oriented subgraph of $_{\preceq}FG_F$ given $S$.

We are interested in sinks and sources of $_{\preceq}S$. 
When confusion may arise, to distinguish relative sources and relative sinks to global sources and global sinks (i.e., to sources and sinks of  $_{\preceq}FG_F$), we say that a source (respectively, sink) of $_{\preceq}S$ is an $S$-source (respectively, $S$-sink).
\begin{lem}
\label{alternato}
Let $F$ be an arbitrary flipclass. 
Each residue $S$ of $F$ has as many $S$-sources as $S$-sinks, and these alternate around the cycle $S$. Furthermore, an $S$-source and an $S$-sink are at distance at least $2$.
\end{lem}
\begin{proof}
The first statement  is true for any oriented cycle. Let us prove the second statement. If $S$ is a non-adjacent residue, then the result is clear. If $S$ is an adjacent residue, then $S$ corresponds to a $3$-flipclass by Proposition~\ref{corrisponde}, and the assertion follows since it is actually true for any $h$-flipclass, with $h>2$. Indeed, let $\Gamma = (x_0 \stackrel{t_1}{\longrightarrow} x_1 \stackrel{t_2}{\longrightarrow}  \cdots \stackrel{t_{h-1}}{\longrightarrow} x_{h-1}\stackrel{t_h}{\longrightarrow}  x_h)$  be a decreasing  path. After performing the $i$-th flip, we have a path that is decreasing everywhere except at position $i$ by Lemma~\ref{espomariancona24}.
\end{proof}

\begin{defn}
Let $F$ be an arbitrary flipclass. Let $\preceq$ be a reflection ordering. We denote by $c_{\preceq}(S)$ the number of $S$-sinks (or $S$-sources) of the residue $S$ of $_{\preceq}FG_F$. In the case that $F$ is of finite type, the number $c_{\preceq}(S)$ does not depend on the chosen reflection ordering, and thus we denote it more simply $c(S)$.
\end{defn}

Adjacent residues correspond to $3$-flipclasses and the two notions of $c$ coincide.
\begin{lem}
\label{corrispondeS}
Let $F$ be an arbitrary flipclass. Let $\preceq$ be a reflection ordering.  Let $S$ be an adjacent residue of  $F$ with associated segment subflipclass $O$. Then 
$$c_{\preceq}(S)=c_{\preceq}(\bar{O}).$$
\end{lem}
\begin{proof}
It is sufficient to observe that the oriented graphs underlying $_{\preceq}S$ and $_{\preceq}FG_{\bar{O}}$ coincide. 
\end{proof}

\begin{lem}
\label{2-3}
Let $F$ be an arbitrary flipclass. Let $\preceq$ be a reflection ordering.
Let $S$ be a residue of $F$ that is a $2m$-cycle with $m\in \{2,3\}$. Then $c_{\preceq}(S)=1$.
\end{lem}
\begin{proof}
Straightforward from Lemma~\ref{alternato}.
\end{proof}

\begin{pro}
\label{2,3,4,5,6,7,9}
Let $S$ be a residue of a flipclass $F$ of finite type.
\begin{enumerate}

\item If $F$ is of dihedral type, then $S$ is a $4$-cycle and $c(S)=1$.

\item If $F$ is of Weyl type,  then $S$ is a $2m$-cycle with
$
m\in\{2,3,4,5,6,7,9\}
$,
and 
$$ c(S)= \left\{ \begin{array}{ll}
1, & \mbox{if $m\in\{2,3,4\}$ }; \\
2, & \mbox{if $m\in\{5,6\}$}; \\
3, & \mbox{if $m\in\{7,9\}$}. 
\end{array} \right. $$

\item If $F$ is of type $H$, then $S$ is a $2m$-cycle with  $2\leq m \leq 13$, and 
$$ c(S)= \left\{ \begin{array}{ll}
1, & \mbox{if $m\in\{2,3,4\}$}; \\
1 \text{ or } 2, & \mbox{if $m =5$}; \\
2, &  \mbox{if $m = 6$; } \\
2 \text{ or } 3, & \mbox{if $m=7$};\\
3, & \mbox{if $m =8$}; \\
3 \text{ or } 4, &  \mbox{if $m \in \{ 9,10\}$}; \\
4 \text{ or } 5, & \mbox{if $m=11$};\\
4, & \mbox{if $m =12$}; \\
4 \text{ or } 5 &  \mbox{if $m = 13$ }. 
\end{array} \right. $$
\end{enumerate}
\end{pro}
\begin{proof}
The first statement follows by  \cite[Section~4]{EMS}.

Let us prove the second and the third statements together.

If $S$ is a $2m$-cycle with $m\in \{2,3\}$, then $c(S)=1$ by Lemma~\ref{2-3}.
If $S$ is a $2m$-cycle with $m\notin \{2,3\}$, then it is an adjacent residue and corresponds to a $3$-flipclass $\bar{O}$ of a Weyl group or of a group of type $H$, by Proposition~\ref{corrisponde}.
The reflection subgroup $W(\bar{O})$ generated by the reflections labeling the paths in $\bar{O}$ is irreducible since, otherwise,  $\bar{O}$ would correspond to a $3$-crown and $S$ would be a $6$-cycle. 
By \cite[Theorem~6.10]{EMS}, we need to check the 3-flipclasses of an irreducible  Coxeter group $W'$ of rank $ \leq 3$ that is either a Weyl group,  or  $H_3$. 
%The reducible groups of rank $\leq 3$ have only flipclasses of cardinality either $4$ or $6$ (they are, indeed,  either $2$- or $3$-crowns). Hence, $W'$ is irreducible.  
An exhaustive check for  $W'$ equal to  $A_3$,   $B_3$, or  $H_3$ gives the result. 
\end{proof}
\begin{rmk}
 If $F$ is of type $A$, then a residue $S$ of $F$ is a $2m$-cycle with
$
m\in\{2,3,4,5\}$.

\end{rmk}

\subsection{Normal Flip Combinatorial Invariance Conjecture}
We recall from the introduction the conjecture driving this work.
\begin{con}
If two flipclasses $F$ and $F'$ of Weyl type  have isomorphic flipgraphs, then $c(F)=c(F')$. 
\end{con}

\begin{rmk}
As can be seen from Proposition~\ref{2,3,4,5,6,7,9},  the Normal Flip Combinatorial Invariance does not hold  for type $H$. Actually, also Flip Combinatorial Invariance fails for type $H$  (see \cite [Remark~7.1]{EMS}).
\end{rmk}

\section{An upper bound}

In this section, we establish an upper bound on the number of sinks of an oriented flip graph in terms of the sinks of its residues.

Let $F$ be an arbitrary flipclass. Fix a reflection ordering $\preceq$. 
\begin{lem}\label{perlemmino1}
Let $s$ be a (global) sink. Then, there exists a residue $S$ and an $S$-sink $p$ such that $\Omega(p)= \{s\}$.
\end{lem}
\begin{proof}
We may clearly suppose that there is more than one global sink.

By Lemma~\ref{bipartito-connesso}, there exists a vertex $v$ such that $|\Omega(v)|>1$ and $s\in \Omega(v)$. Let $v$ be a minimal such vertex with respect to the partial order of Lemma~\ref{loreto} where the sinks are the (globally) minimal elements. %Explicitly, this means that $|\Omega(v)|>1$, the sink $s$ is downstream from $v$ and, for any vertex $v'$ strictly downstream from $v$, either $s\notin  \Omega(v')$ or $\{s\} =  \Omega(v')$.
Let $\alpha$ and $\beta$ be two oriented paths starting at $v$ such that $\alpha$ ends in $s$ and $\beta$ ends in a different sink. Let $v \stackrel{i}{\rightarrow} a$ and $v \stackrel{j}\rightarrow b$ be the first steps respectively in the paths $\alpha$ and $\beta$. Consider the $\{i,j\}$-residue $S$ containing $v$. Note that $v$ is an $S$-source. There are two $S$-sinks that are downstream from $v$. One, call it $q$,  is downstream from $a$ and hence $\Omega(q)\subseteq \Omega(a)$. We can choose $p=q$ since $ \Omega(a)=\{s\}$.
\end{proof}

%\begin{cor}
%Let $G$ be the graph with $\{S: \text{S residue}\}$ as vertex set and where $S $ and $S'$ are adjacent if $\tilde{\Omega}(S) \cap  \tilde{\Omega}(S')\neq \emptyset$. Then $G$ is connected
%\end{cor}
%\begin{proof}
%Towards a contradiction, suppose that $G$ is not connected. Let $A$ and $B$ be two sets of residues such that $A$ Then there exist two sets of global sinks $\mathcal S_1$ and $\mathcal S_2$
%\end{proof}

For the rest of this section, fix a map $\tilde{\Omega}$ that sends a vertex $v$ of $_{\preceq}FG_F$ to a (unique, global) sink  $\tilde{\Omega}(v)$ in the set $\Omega(v)$.
Note that, in general, $|\Omega(v)|\geq 1$, so $\tilde{\Omega}$ depends on the choice we make. 

For each residue $S$, we let
$$\tilde{\Omega}(S) =\{ \tilde{\Omega}(v) : \text{$v$ is an $S$-sink} \}.$$

Recall that $\mathcal S_{u}$ (respectively, $\mathcal S_{d}$) denotes the set of sources (respectively, sinks) of $_{\preceq} FG_F$. Recall also $|\mathcal S_{d}|=c_{\preceq}(F)$. 

\begin{lem}\label{perlemmino2}
Let $U$ be a proper nonempty subset of $\mathcal S_{d}$. Then there exists a residue $S$ such that $\tilde{\Omega}(S)\cap U \neq \emptyset$ and $ \tilde{\Omega}(S) \cap (\mathcal S_{d}\setminus U) \neq \emptyset$.
\end{lem}
\begin{proof}
By Lemma~\ref{bipartito-connesso}, there exists a vertex $v$ such that $\Omega(v)\cap U \neq \emptyset$ and $\Omega(v)\cap (\mathcal S_{d}\setminus U) \neq \emptyset$. 
Let $v$ be a minimal such vertex. 
Let $\alpha$ and $\beta$ be two oriented paths starting at $v$ such that $\alpha$ ends in a sink of $U$ and $\beta$ ends in a sink in $(\mathcal S_{d}\setminus U)$.  Let $v \stackrel{i}{\rightarrow} a$ and $v \stackrel{j}\rightarrow b$ be the first steps respectively in the paths $\alpha$ and $\beta$. Consider the $\{i,j\}$-residue $S$ containing $v$. Note that $v$ is an $S$-source. There are two $S$-sinks that are downstream from $v$. One, call it $p$,  is downstream from $a$ and the other, call it $q$, is downstream from $b$. Then $\Omega(p)\subseteq \Omega(a)$ and  $\Omega(q)\subseteq \Omega(b)$. 
By minimality, $\Omega(a) \subseteq U$ and $\Omega(b) \subseteq (\mathcal S_{d}\setminus U)$. Hence  $\tilde{\Omega}(p) \in U$ and $\tilde{\Omega}(q) \in (\mathcal S_{d}\setminus U)$. Hence the residue $S$ satisfies the required properties.
\end{proof}

We shall apply the following elementary lemma.
\begin{lem}
\label{lemmino}
Let $X$ be a finite set and let $\{X_{\alpha} : \alpha\in A\}$ be a finite family of nonempty subsets of $X$ such that 
\begin{enumerate}
\item
\label{lemmino1}
 $\cup_{\alpha \in A} X_{\alpha} = X$ and, 
\item 
\label{lemmino2}
for all proper nonnempty subset $U$ of  $X$, there exists $\alpha \in A$ such that $X_{\alpha} \cap U \neq \emptyset$ and    $X_{\alpha} \cap (X\setminus U) \neq \emptyset$.
\end{enumerate}
Then $$|X|\leq 1+ \sum_{\alpha \in A} (|X_{\alpha}|-1).$$
\end{lem}
\begin{proof}
Order the elements in the following way. Let $A=\{\alpha_1, \ldots, \alpha_{|A|}\}$ be such that, for each $i\in \{1,\ldots,|A|-1\}$, the set  $X_{\alpha_{i+1}}$ satisfies $X_{\alpha_{i+1}} \cap U_i \neq \emptyset$ 
%and    $X_{\alpha_{i+1}} \cap (X\setminus U_i) \neq \emptyset$
 for $U_i=\cup_{k=1}^{i} X_{\alpha_k}$. Such an order exists by (\ref{lemmino2}).
Then 
$$|X| =  |X_{\alpha_1}| + \sum_{i=2}^{|A|}  |X_{\alpha_i}\setminus U_{i-1}| \leq 1+ \sum_{\alpha \in A} ( |X_{\alpha}|-1 ) ,$$
as desired.
\end{proof}

\begin{thm}
\label{minore1}
Let $F$ be an arbitrary flipclass of an arbitrary Coxeter group $W$. Fix a reflection ordering $\preceq$ of $W$. Then
$$c_{\preceq}(F) \leq 1 + \sum_{S}(c_{\preceq}(S)-1),$$
where the sum is over the residues of the flipgraph of $F$.
\end{thm}

\begin{proof}
By Lemmas~\ref{perlemmino1} and~\ref{perlemmino2}, we can apply Lemma~\ref{lemmino} to the nonempty sets $\tilde{\Omega}(S)$, which satisfy $|\tilde{\Omega}(S)|\leq c(S)$.
\end{proof}

\begin{thm}
\label{minore2}
Let $F$ be a flipclass. 
\begin{enumerate}
\item
\label{minoreft}
 If $F$ is of finite type, then 
$$c(F) \leq 1 + \sum_{S}(c(S)-1),$$
where the sum is over the residues of the flipgraph of $F$.
\item
\label{minorew}  If $F$ is of Weyl type, then 
$$c(F) \leq 1 + \sum_{S}(c(S)-1),$$
where the sum is over the residues of the flipgraph of $F$ and
$$ c(S)= \left\{ \begin{array}{ll}
1, & \mbox{if $S$ is a $2m$-cycle with $m\in\{2,3,4\}$; } \\
2, & \mbox{if $S$ is a $2m$-cycle with $m\in\{5,6\}$; } \\
3, & \mbox{if $S$ is a $2m$-cycle with $m\in\{7,9\}$. } 
\end{array} \right. $$
\end{enumerate}
\end{thm}

\begin{proof}
The first statement is just Theorem~\ref{minore1} for flipclasses of finite type. The second statement follows from the first and Proposition~\ref{2,3,4,5,6,7,9}.
\end{proof}

\begin{rmk}
Theorem~\ref{minore1} holds for arbitrary flipclasses. In this generality, we do not know whether the number of increasing paths of $F$ is independent of the chosen reflection ordering. Theorem~\ref{minore2}, (\ref{minoreft}), holds for flipclasses of finite type, for which we know that  the number of increasing paths of $F$ is independent of the chosen reflection ordering but we also know that the flip combinatorial invariance does not hold (for type $H$). 
Theorem~\ref{minore2}, (\ref{minorew}), holds for flipclasses of Weyl type, for which we know that  the number of increasing paths of $F$ is independent of the chosen reflection ordering and we conjecture  that the flip combinatorial invariance  holds (we know that it holds for $h$-flipclasses up to $h=6$, see \cite{EMS}). 
\end{rmk}

\begin{rmk}
Note that the sum in Theorem~\ref{minore2}, (\ref{minorew}), can be taken over all residues that are $2m$-cycles, with $m\in\{5,6,7,9\}$. These are adjacent residues and far fewer than the residues that are $2m$-cycles, with $m\in\{2,3,4\}$. 
\end{rmk}

\begin{rmk}
The fact that residues play a fundamental role mirrors the fact that (parabolic) dihedral subgroups and dihedral intervals play a fundamental role in understanding the combinatorial properties of Coxeter groups (see, for example, \cite{BCM1}, \cite{CMEur}, \cite{CMDis}, \cite{DyeComp1}, \cite{DyeComp}, \cite{DyePre}, \cite{Mtrans}, \cite{M}).
\end{rmk}

\section{A lower bound}

This short section contains a lower bound on the number of increasing paths of a flipclass.

\begin{thm}\label{terni}
Let $F$ be an $h$-flipclass of finite type from $u$ to $v$. Let $u< x_1< x_2 < \cdots < x_r<  v$ and  $0< i_1 < i_2 < \cdots < i_r < h$.  Let $J=[h-1]\setminus \{i_k : k\in [r]\}$. Then
$$c(F) \geq \sum_{O}  c(O),$$
where $O$ runs over the  $J$-subflipclasses of $F$ of paths passing through $x_k$ at time $i_k$, for $k\in[r]$. 
\end{thm}
\begin{proof}
Suppose $r=1$, and let $x=x_1$. By  \cite[Proposition~6.7]{EMS}, we can choose any reflection ordering. Let us choose a reflection ordering $\preceq_x$ having $ \{t \in T : \ell(tx) < \ell(x)\}$ as an initial segment, which exists by \cite[Lemma~2.11]{DyeComp}.  Then every  increasing path ending at $x$ concatenates increasingly with every increasing path starting at $x$, and the assertion is proved.
If $r>1$, we first apply the previous case with $x=x_1$ to obtain $c(F) \geq  \sum_{O^{x_1}}  c(O^{x_1})$. Each $O^{x_1}$ in the sum satisfies  $O^{x_1}=\bar{O}^{x_1}_1 \times \bar{O}^{x_1}_2$. Then we can apply the same argument  to $\bar{O}^{x_1}_2$ for $x=x_2$ with $\preceq_{x_2}$ (being $F$ of finite type, every choice of the ordering is fine). Iterating the process, we get the assertion. 
\end{proof}

\begin{rmk}
The reflection ordering appearing in the proof of Theorem~\ref{terni} is called \emph{Deodhar with respect to $x$} in \cite[Definition~3.3]{BGL}. An argument analogous to that of the proof of Theorem~\ref{terni} has been used when showing that shortcuts of particular hypercube decompositions are calculating (see \cite{B-G}, \cite{B-M}). The idea of iterating the basic step in that framework led to the notion of doubleshortcut (see \cite{EM2}). 
\end{rmk}

\begin{rmk}
\label{cracker}
Theorem~\ref{terni} implies that, for a segment subflipclass $O$ of a flipclass $F$ of finite type, we have $c(\bar{O})\leq c(F)$. This is a monotonicity result on the number of increasing paths. (Observe, however, that  an increasing path of $\bar{O}$ in general cannot be extended to an increasing path of $F$). In particular, given a residue $S$, we have $c(S)\leq c(F)$ by Proposition~\ref{corrisponde} and  Lemmas~\ref{corrispondeS} and~\ref{2-3}. 
\end{rmk}

\section{Consequences}
In this section, we collect the main results stemming from the framework developed in the previous sections. 
\begin{cor}
\label{inug}
Let $F$ be an $h$-flipclass of finite type. Then
$$\max_{(x,i)} \left( \sum_{O}  c(O) \right) \leq c(F) \leq 1 + \sum_{S}(c(S)-1),$$
where the first sum is over all maximal subflipclasses of $F$ of paths passing through $x$ at time $i$  and the second sum  is over the residues of the flipgraph of $F$. 
\end{cor}
\begin{proof}
Straightforward from Theorem~\ref{minore2}, (\ref{minoreft}) and Theorem~\ref{terni}.
\end{proof}

\begin{cor}
Every flipclass $F$ of finite type with $c(F)=1$ is  normal.
\end{cor}
\begin{proof}
Since $c(O) \geq 1$%by Proposition~\ref{greedy}
, the equality $c(F) = 1$ implies that the sum in Theorem~\ref{terni} is over a singleton.
\end{proof}

\begin{thm}
\label{c=1}
Let $F$ be a flipclass of finite type. Then  $c(F)=1$ if and only if each residue $S$ of $F$ satisfies $c(S)=1$. In particular, if each residue of $F$ is a $2m$-cycle with $m\leq 4$, then $c(F)=1$.
\end{thm}
\begin{proof}
Let us prove the first assertion. If $c(F)=1$, then $c(S)=1$ for each residue  $S$ by Remark~\ref{cracker}. The converse is straightforward from Corollary~\ref{inug}. 

The second assertion follows by Proposition~\ref{2,3,4,5,6,7,9}.
\end{proof}
\begin{rmk}
Let $F$ be an arbitrary flipclass (i.e., not necessarily of finite type) and let $\preceq$ be a reflection ordering. By Theorem~\ref{minore1}, we have that, if $c_{\preceq} (S)=1$ for each residue $S$, then $c_{\preceq}(F)=1$, since every flipclass has at least one increasing path with respect to any reflection ordering. 
\end{rmk}

The sufficient condition of Theorem~\ref{c=1} becomes also necessary for flipclasses of Weyl type. This implies that  the property ``$c(F)=1$'' is determined by the flipgraph. Hence, it provides substantial evidence in favour of Normal Flip Combinatorial Invariance. 
\begin{cor}
Let $F$ be a flipclass of Weyl type. Then  $c(F)=1$ if and only if each residue of $F$ is a $2m$-cycle, with $m\leq 4$.
\end{cor}
\begin{proof}
Follows by Proposition~\ref{2,3,4,5,6,7,9} and Theorem~\ref{c=1}.
\end{proof}
In particular, we have the following result.
\begin{cor}\label{tardi}
Let $F$ and $F'$ be two flipclasses of Weyl type with isomorphic flipgraphs. Then 
$$c(F)=1 \iff c(F')=1.$$
\end{cor}

\begin{rmk}
In fact, the hypotheses in Corollary~\ref{tardi} can be relaxed. It is sufficient that there is an edge-labelled graph isomorphism between the flipgraphs of $F$ and $F'$, i.e. a graph isomorphism inducing an arbitrary permutation on the labels (this happens, in particular, if $F$ and $F'$ are anti-isomorphic, see \cite[Definition~3.8]{EMS}).
\end{rmk}

\begin{figure}[ht]
\centering
\begin{tikzpicture}[scale=0.80,
  v/.style={circle,fill=black,inner sep=0pt,minimum size=3.2pt},
  f1/.style={red!80!black,thick},
  f2/.style={blue!80!black,thick,dashed},
  f3/.style={green!50!black,thick,densely dotted},
  f4/.style={orange!90!black,thick,dash dot}]
  \draw[f4] (6.294,-1.107) -- (5.452,-2.831);
  \draw[f3] (6.294,-1.107) -- (3.586,-1.207);
  \draw[f2] (6.294,-1.107) -- (7.002,1.749);
  \draw[f1] (6.294,-1.107) -- (4.635,-0.336);
  \draw[f3] (5.452,-2.831) -- (3.066,-4.181);
  \draw[f2] (5.452,-2.831) -- (6.145,0.186);
  \draw[f1] (5.452,-2.831) -- (3.869,-2.429);
  \draw[f4] (3.586,-1.207) -- (2.013,-3.288);
  \draw[f2] (3.586,-1.207) -- (2.442,0.953);
  \draw[f1] (3.586,-1.207) -- (2.071,-0.754);
  \draw[f3] (2.013,-3.288) -- (1.404,-5.497);
  \draw[f2] (2.013,-3.288) -- (0.574,-1.262);
  \draw[f1] (2.013,-3.288) -- (-0.216,-2.263);
  \draw[f4] (1.404,-5.497) -- (3.066,-4.181);
  \draw[f2] (1.404,-5.497) -- (-0.126,-3.728);
  \draw[f1] (1.404,-5.497) -- (-0.850,-5.166);
  \draw[f2] (3.066,-4.181) -- (1.517,-1.782);
  \draw[f1] (3.066,-4.181) -- (1.149,-4.414);
  \draw[f4] (7.002,1.749) -- (6.145,0.186);
  \draw[f3] (7.002,1.749) -- (5.087,3.982);
  \draw[f1] (7.002,1.749) -- (5.463,2.561);
  \draw[f3] (6.145,0.186) -- (4.068,2.326);
  \draw[f1] (6.145,0.186) -- (4.197,0.903);
  \draw[f4] (5.087,3.982) -- (4.068,2.326);
  \draw[f2] (5.087,3.982) -- (2.864,3.468);
  \draw[f1] (5.087,3.982) -- (3.438,4.947);
  \draw[f2] (4.068,2.326) -- (1.457,0.917);
  \draw[f1] (4.068,2.326) -- (2.089,3.267);
  \draw[f4] (2.442,0.953) -- (0.574,-1.262);
  \draw[f3] (2.442,0.953) -- (2.864,3.468);
  \draw[f1] (2.442,0.953) -- (-0.193,2.567);
  \draw[f3] (0.574,-1.262) -- (-0.126,-3.728);
  \draw[f1] (0.574,-1.262) -- (-1.761,-0.063);
  \draw[f4] (2.864,3.468) -- (1.457,0.917);
  \draw[f1] (2.864,3.468) -- (0.237,4.195);
  \draw[f3] (1.457,0.917) -- (1.517,-1.782);
  \draw[f1] (1.457,0.917) -- (-1.289,1.785);
  \draw[f4] (-0.126,-3.728) -- (1.517,-1.782);
  \draw[f1] (-0.126,-3.728) -- (-2.528,-2.651);
  \draw[f1] (1.517,-1.782) -- (-1.421,-1.667);
  \draw[f4] (4.635,-0.336) -- (3.869,-2.429);
  \draw[f3] (4.635,-0.336) -- (2.071,-0.754);
  \draw[f2] (4.635,-0.336) -- (5.463,2.561);
  \draw[f3] (3.869,-2.429) -- (1.149,-4.414);
  \draw[f2] (3.869,-2.429) -- (4.197,0.903);
  \draw[f4] (2.071,-0.754) -- (-0.216,-2.263);
  \draw[f2] (2.071,-0.754) -- (0.014,1.793);
  \draw[f3] (-0.216,-2.263) -- (-0.850,-5.166);
  \draw[f2] (-0.216,-2.263) -- (-2.472,-0.101);
  \draw[f4] (-0.850,-5.166) -- (1.149,-4.414);
  \draw[f2] (-0.850,-5.166) -- (-3.230,-4.881);
  \draw[f2] (1.149,-4.414) -- (-2.242,-3.810);
  \draw[f4] (5.463,2.561) -- (4.197,0.903);
  \draw[f3] (5.463,2.561) -- (3.438,4.947);
  \draw[f3] (4.197,0.903) -- (2.089,3.267);
  \draw[f4] (3.438,4.947) -- (2.089,3.267);
  \draw[f2] (3.438,4.947) -- (0.898,5.104);
  \draw[f2] (2.089,3.267) -- (-1.024,3.848);
  \draw[f4] (0.014,1.793) -- (-2.472,-0.101);
  \draw[f3] (0.014,1.793) -- (0.898,5.104);
  \draw[f1] (0.014,1.793) -- (-2.203,3.597);
  \draw[f3] (-2.472,-0.101) -- (-4.803,-0.065);
  \draw[f1] (-2.472,-0.101) -- (-3.654,1.365);
  \draw[f4] (0.898,5.104) -- (-1.024,3.848);
  \draw[f1] (0.898,5.104) -- (-1.461,5.593);
  \draw[f3] (-1.024,3.848) -- (-3.991,2.022);
  \draw[f1] (-1.024,3.848) -- (-3.139,4.132);
  \draw[f4] (-4.803,-0.065) -- (-3.991,2.022);
  \draw[f2] (-4.803,-0.065) -- (-5.460,-2.779);
  \draw[f1] (-4.803,-0.065) -- (-6.294,1.048);
  \draw[f2] (-3.991,2.022) -- (-4.322,-1.261);
  \draw[f1] (-3.991,2.022) -- (-5.492,2.715);
  \draw[f4] (-3.230,-4.881) -- (-2.242,-3.810);
  \draw[f3] (-3.230,-4.881) -- (-5.460,-2.779);
  \draw[f1] (-3.230,-4.881) -- (-4.874,-3.881);
  \draw[f3] (-2.242,-3.810) -- (-4.322,-1.261);
  \draw[f1] (-2.242,-3.810) -- (-3.955,-2.377);
  \draw[f4] (-5.460,-2.779) -- (-4.322,-1.261);
  \draw[f1] (-5.460,-2.779) -- (-7.012,-1.866);
  \draw[f1] (-4.322,-1.261) -- (-6.150,-0.307);
  \draw[f4] (-0.193,2.567) -- (-1.761,-0.063);
  \draw[f3] (-0.193,2.567) -- (0.237,4.195);
  \draw[f2] (-0.193,2.567) -- (-2.203,3.597);
  \draw[f3] (-1.761,-0.063) -- (-2.528,-2.651);
  \draw[f2] (-1.761,-0.063) -- (-3.654,1.365);
  \draw[f4] (0.237,4.195) -- (-1.289,1.785);
  \draw[f2] (0.237,4.195) -- (-1.461,5.593);
  \draw[f3] (-1.289,1.785) -- (-1.421,-1.667);
  \draw[f2] (-1.289,1.785) -- (-3.139,4.132);
  \draw[f4] (-2.528,-2.651) -- (-1.421,-1.667);
  \draw[f2] (-2.528,-2.651) -- (-4.874,-3.881);
  \draw[f2] (-1.421,-1.667) -- (-3.955,-2.377);
  \draw[f4] (-2.203,3.597) -- (-3.654,1.365);
  \draw[f3] (-2.203,3.597) -- (-1.461,5.593);
  \draw[f3] (-3.654,1.365) -- (-6.294,1.048);
  \draw[f4] (-1.461,5.593) -- (-3.139,4.132);
  \draw[f3] (-3.139,4.132) -- (-5.492,2.715);
  \draw[f4] (-6.294,1.048) -- (-5.492,2.715);
  \draw[f2] (-6.294,1.048) -- (-7.012,-1.866);
  \draw[f2] (-5.492,2.715) -- (-6.150,-0.307);
  \draw[f4] (-4.874,-3.881) -- (-3.955,-2.377);
  \draw[f3] (-4.874,-3.881) -- (-7.012,-1.866);
  \draw[f3] (-3.955,-2.377) -- (-6.150,-0.307);
  \draw[f4] (-7.012,-1.866) -- (-6.150,-0.307);
  \node[v] (n1) at (6.294,-1.107) {};
  \node[v] (n2) at (5.452,-2.831) {};
  \node[v] (n3) at (3.586,-1.207) {};
  \node[v] (n4) at (2.013,-3.288) {};
  \node[v] (n5) at (1.404,-5.497) {};
  \node[v] (n6) at (3.066,-4.181) {};
  \node[v] (n7) at (7.002,1.749) {};
  \node[v] (n8) at (6.145,0.186) {};
  \node[v] (n9) at (5.087,3.982) {};
  \node[v] (n10) at (4.068,2.326) {};
  \node[v] (n11) at (2.442,0.953) {};
  \node[v] (n12) at (0.574,-1.262) {};
  \node[v] (n13) at (2.864,3.468) {};
  \node[v] (n14) at (1.457,0.917) {};
  \node[v] (n15) at (-0.126,-3.728) {};
  \node[v] (n16) at (1.517,-1.782) {};
  \node[v] (n17) at (4.635,-0.336) {};
  \node[v] (n18) at (3.869,-2.429) {};
  \node[v] (n19) at (2.071,-0.754) {};
  \node[v] (n20) at (-0.216,-2.263) {};
  \node[v] (n21) at (-0.850,-5.166) {};
  \node[v] (n22) at (1.149,-4.414) {};
  \node[v] (n23) at (5.463,2.561) {};
  \node[v] (n24) at (4.197,0.903) {};
  \node[v] (n25) at (3.438,4.947) {};
  \node[v] (n26) at (2.089,3.267) {};
  \node[v] (n27) at (0.014,1.793) {};
  \node[v] (n28) at (-2.472,-0.101) {};
  \node[v] (n29) at (0.898,5.104) {};
  \node[v] (n30) at (-1.024,3.848) {};
  \node[v] (n31) at (-4.803,-0.065) {};
  \node[v] (n32) at (-3.991,2.022) {};
  \node[v] (n33) at (-3.230,-4.881) {};
  \node[v] (n34) at (-2.242,-3.810) {};
  \node[v] (n35) at (-5.460,-2.779) {};
  \node[v] (n36) at (-4.322,-1.261) {};
  \node[v] (n37) at (-0.193,2.567) {};
  \node[v] (n38) at (-1.761,-0.063) {};
  \node[v] (n39) at (0.237,4.195) {};
  \node[v] (n40) at (-1.289,1.785) {};
  \node[v] (n41) at (-2.528,-2.651) {};
  \node[v] (n42) at (-1.421,-1.667) {};
  \node[v] (n43) at (-2.203,3.597) {};
  \node[v] (n44) at (-3.654,1.365) {};
  \node[v] (n45) at (-1.461,5.593) {};
  \node[v] (n46) at (-3.139,4.132) {};
  \node[v] (n47) at (-6.294,1.048) {};
  \node[v] (n48) at (-5.492,2.715) {};
  \node[v] (n49) at (-4.874,-3.881) {};
  \node[v] (n50) at (-3.955,-2.377) {};
  \node[v] (n51) at (-7.012,-1.866) {};
  \node[v] (n52) at (-6.150,-0.307) {};
\end{tikzpicture}
\caption{A flipgraph of a 5-flipclass with $c(F) = 1$}
\label{arriva}
\end{figure}
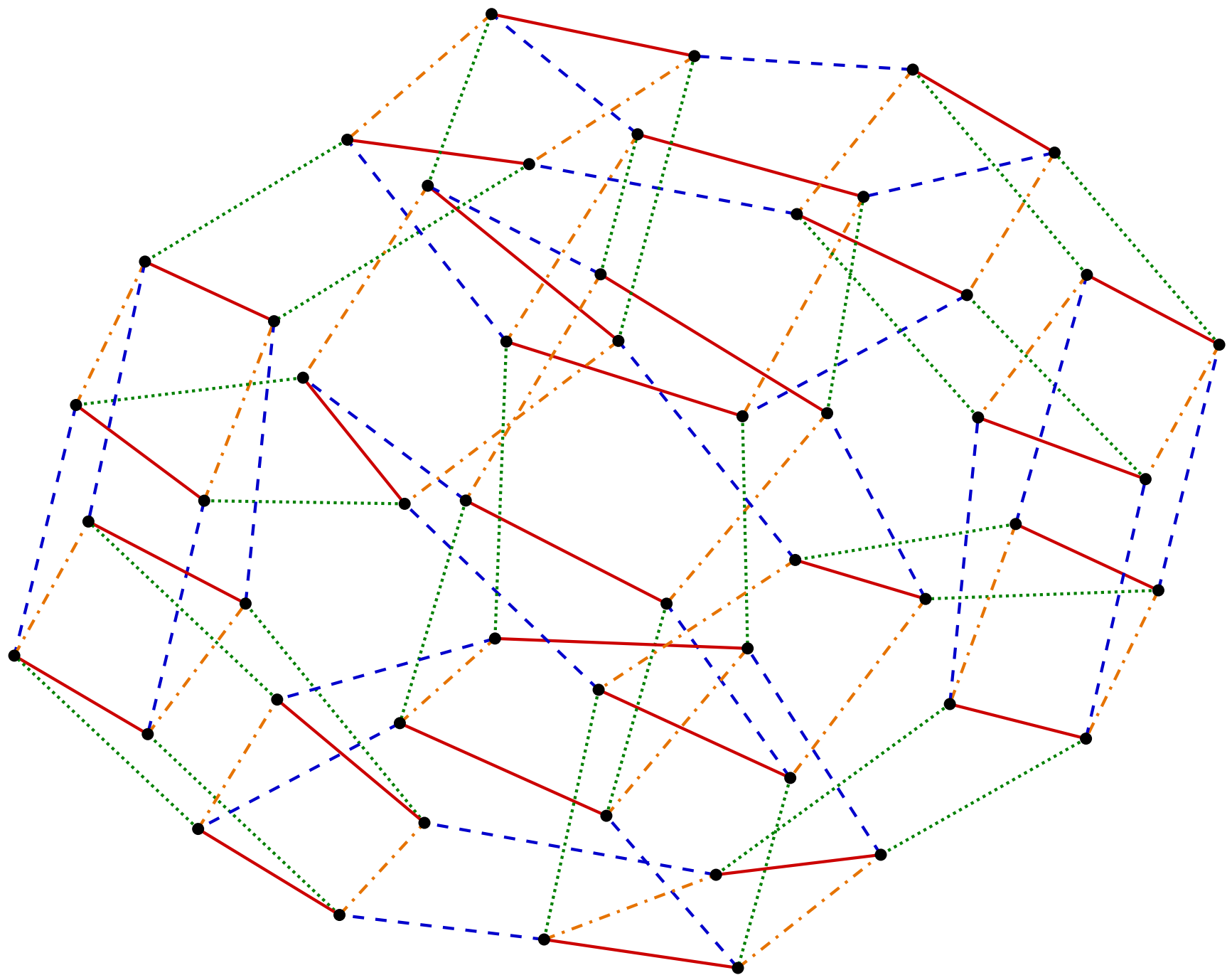

\begin{exa}
Let $F$ be a flipclass whose flipgraph is isomorphic (as an edge-labelled graph) to the one depicted in Figure~\ref{arriva}. Since all residues are $2m$-cycles with $m \leq 4$, then $c(F) = 1$. An example of such an $F$ is the flipclass of the path 
$ 12354 \longrightarrow 14352 \longrightarrow 15342 \longrightarrow 25341$ in $S_5$.
\end{exa}

\bigskip
{\bf Acknowledgments:} 
The authors are members of the INDAM group GNSAGA.

\end{document}